\documentclass{article}

\usepackage{amsmath}
\usepackage[francais]{babel}
\usepackage{amssymb}
\usepackage[applemac]{inputenc}

\def\F{\bf F}
\def\C{\bf C}
\def\R{\bf R}
\def\G{\mathcal{G}}
\def\Weyl{\rm Weyl}
\def\Sym{\rm Sym}

\def\Inv{\rm Inv}
\def\n{\noindent}
\def\>{\rangle}
\def\<{\langle}

\begin{document}
 
 \begin{center} 
  
  Cohomological invariants mod 2 of Weyl groups
  
  \medskip
  Jean-Pierre Serre
  \end{center}

  \bigskip
 
 Let $G$ be the Weyl group of a root system, i.e., a crystallographic finite Coxeter group, cf. [B], chap.VI, §4.1. Let $k_0$ be a field of characteristic $≠2$, let $H^\bullet(k_0) = \oplus_{n \geqslant 0} \ H^n(k_0,{\bf F}_2)$ and let $I_G = {\rm Inv}_{k_0}(G)$ be the ring of cohomological invariants mod 2 of $G$, as defined in [S], §4; it is a graded $H^\bullet(k_0)$-algebra. When $G$ is of
 type $\sf{A}$, it is isomorphic to a symmetric group $\Sym$$_n$, and  $I_G$ is $H^\bullet(k_0)$-free of rank $1+[n/2]$, with an explicit basis $w_0=1, w_1,...,w_{[n/2]}$, cf. [S], chap.VII.

 \smallskip
   In order to extend this description of $I_G$ to the general case, define $S_G$ to be the set of elements $g\in G$ with $g^2=1$; an element of $S_G$ shall be called an {\it involution} of $G$. Let $\Sigma_G$ be the set of conjugation classes of elements
   of $S_G$. 
   
   \smallskip
   \n {\bf Theorem A.} {\it There exists a natural injection $e :\Sigma_G \to I_G$ whose image is an $H^\bullet(k_0)$-basis of $I_G$.}
   
   [Equivalently : the $H^\bullet(k_0)$-module $I_G$ is canonically isomorphic to the set of all maps $\Sigma_G \to H^\bullet(k_0)$.]
   
   \smallskip
     The map $e$ is compatible with the grading of $I_G$: if $g\in S_G$,
     define the {\it degree} of $g$ to be be the multiplicity of $-1$ as an eigenvalue of $g$ in the standard  linear representation of $G$ as a Coxeter group; let $\Sigma_{G,n}$ be the set of involution classes of degree $n$. 
    {\it If $\sigma \in \Sigma_{G,n}$, then $e(\sigma)$ belongs to the $n$-th component $I^n_G$ of $I_G$}.
     
\medskip

   \n {\it Examples.} 1. When $G = \Sym$$_n$, the elements of $\Sigma_G$ are the conjugation classes of the products of $i$ disjoint transpositions, with $2i \leqslant n$, and we recover the fact that $H^\bullet(k_0)$-free of rank $1+[n/2]$, with a basis made up of elements of degree
  $ 0,1,..., [n/2]$. In that case the canonical basis is made up of the $w_i^{\rm gal}$, which are closely related to the $w_i$ mentioned above, cf. [S], §25.\\
   2. When $G = \Weyl(\sf{E}_8)$, we have $|\Sigma_{G,n}|=1$ for $0 \leqslant n \leqslant8$, with the only exception of
   $n=4$ where $|\Sigma_{G,n}|=2$; and, of course, $\Sigma_{G,n}=\varnothing$ for $n>8$. Hence $I_G$ is a free $H^\bullet(k_0)$-module of rank 10, with a basis made up of elements of degree $0,1,2,3,4,4,5,6,7,8$.\\  
   3. For $\sf{E}_7$ and $\sf{E}_6$, the degrees are  $0,1,2,3,3,4,4,5,6,7$ and $0,1,2,3,4$.
   
   \bigskip
   \n {\bf Definition of the map} $e :\Sigma_G \to I_G$.
   
   \smallskip
   
     Let $a$ be an element of $I_G$ and let $g$ be an involution of $G$ of degree $n$. We first define a « scalar product »
  $  \< a,g \>,  $ which is an element of $H^\bullet(k_0)$. To do so, choose a splitting 
  $g = s_1\cdots s_n$, where the $s_i$ are commuting reflections (recall that a reflection is an involution of degree 1);
  such a splitting always exists. Let $C = \< s_1,...,s_n\>$ be the group generated by the $s_i$, and let $a_C \in I_C$ be the 
  image of $a$ by the restriction map $I_G \to I_C$. The algebra $I_C$ has a natural basis $(\alpha_I)$ indexed by the subsets $I$ of $[1,n]$, cf. [S], §16.4.
  Let $\lambda_C \in H^\bullet(k_0)$ be the coefficient of $\alpha_{[1,n]}$ in $a_C$ ( « top coefficient »). One can show that $\lambda_C$ {\it is independent of the chosen splitting of $g$, i.e., that it only depends on $a$ and} $g$. We then define the scalar product $\< a,g \> $ as $\lambda_C$; we have $\< a,g \> = \< a,g' \>$ if $g$ and $g'$ are conjugate in $G$; this allows us to define $\< a,\sigma \>$ for every $\sigma \in \Sigma_G$. For a given $\sigma$, the map $a \mapsto \< a,\sigma \>$ is $H^\bullet(k_0)$-linear; if $a$ has degree $m$, then $\< a,\sigma \>$ has degree $m-n$ (one may view $a \mapsto \< a,\sigma \>$ as an $n$-th fold residue map).
  
\n{\it Example.} Choose for $a$ a Stiefel-Whitney  class $w_i^{\rm gal}(\rm Cox)$ of the Coxeter representation of $G$.
One has $\<a,\sigma\> = 0$ if $i \neq \deg(\sigma)$ and $\<a,g\> = 1$ if $i = \deg(\sigma)$.
  
  \medskip
  
  \n {\bf Theorem B.} 
   
  (i) {\it If $a\in I_G$ is such that $\< a,\sigma \> = 0$ for every $\sigma$, then $a=0$.}
  
  (ii) {\it Let $n$ be an integer. For every $\sigma \in \Sigma_G$ of degree $n$, there exists $e(\sigma) \in I^n_G$
            such that $\< e(\sigma)
            ,\sigma \> = ~1$ and $\< e(\sigma),\sigma' \> = 0$ for every $\sigma' \neq \sigma$.}
            
        \n   [Note that, by (i), such an $e(\sigma)$ is unique.]
          
          \medskip
          
          It is clear that Theorem B implies Theorem A.
          
    \medskip
    
    \n {\it Indications on the proof of part} (i){\it of Theorem B}.
    
    \smallskip
    An induction argument shows that, if $\< a,\sigma \> = 0$ for every $\sigma$, then the restriction of $a$ to every « cube » (i.e., subgroup 
    generated by commuting reflections) is $0$. In that case, if the characteristic of $k_0$ is good for $G$, the arguments of [S], §25, show that $a=0$. This already covers the case where the irreducible components of $G$ are of classical type, since every characteristic $\neq 2$ is good. The exceptional
    types can be reduced to the classical ones, thanks to the fact that, if $G$ is such a Weyl group, there exists a subgroup $G'$ of $G$, generated by a subset of $S_G$ (hence also a Weyl group), which is of classical type, and has {\it odd index} in $G$: for  $G$ of type $\sf{E}_6,\sf{E}_7, \sf{E}_8, \sf{F}_4, \sf{G}_2$, one takes $G'$ of type
 $ \sf{D}_5, \sf{A}_1\times \sf{D}_6, \sf{D}_8, \sf{B_4}, \sf{A}_1 \times \sf{A}_1$, respectively; one has $(G:G') = 27, 63, 135, 3,  3.$ One then uses the fact that the restriction map $I_G \to I_{G'}$ is injective, cf. [S], prop.14.4, and that every cube of 
    $G$ is conjugate to a cube of $G'$.
    
    \medskip
         
           \n {\it Indications on the proof of part} (ii){\it of Theorem B}.
    
    \smallskip
  We need to construct enough cohomological invariants. For most Weyl groups, this is done by using Stiefel-Whitney classes. 
  For instance, for Weyl$(\sf{E}_6)$, one takes the $w_i^{\rm gal}(\rm Cox)$, $i = 0,1,2,3,4$. There are however
 three cases where we have to do otherwise. For each one, there are two distinct classes of involutions $\sigma, \sigma'$ of the same degree $n$ for which it is hard to find $a \in I^n_G$ with $\< a,\sigma \> = 0, \ \< a,\sigma' \> =1$. These cases are :
  $\sf{D}_{2{\it n}}$, $n \geqslant 3$; $\sf{E}_7$, $n = 3$ and $4$; $\sf{E}_8$, $n=4$.  
  
  \smallskip
  
  \n  For those, we use the relation given by Milnor's conjecture (now Voevodsky's theorem) between Witt invariants and cohomological invariants mod 2. The method applies to every linear
    group $\G$ over $k_0$. The ring $\Inv$$_{k_0}(\G,$$W)$ of Witt invariants of $\G$ (as defined in [S], §27.3) has a natural filtration : an invariant $h$ has filtration $\geqslant n$ if, for every extension $k/k_0$ and every $\G$-torsor $t$ of $\G$ over  $k$, the element $h(t)$ of the Witt ring $W(k)$ belongs to the $n$-th power of the canonical ideal of $W(k)$; in that case, $h$ defines (via the Milnor construction) an element $a_h$ of ${\rm Inv}$$_{k_0}^n(\G,{\bf F}_2)$ which is $0$ if and only if the filtration of $h$ is $>n$. We thus get {\it an injective map} \  $ {\rm gr}^n \Inv$$_{k_0}(\G,$$W) \ \to$ \ ${\rm Inv}$$_{k_0}^n(\G,{\bf F}_2)$.
             
            \smallskip
            
            We apply this to $\G = G$, where $G$ is as in the three cases above. One can find a linear orthogonal representation
            of $G$ whose Brauer character  $\chi$ is such that $\chi(\sigma) - \chi(\sigma') = 2^n$. This gives a $G$-quadratic form,
            hence an element of $\Inv$$_{k_0}(G,W)$; one modifies slightly that element to make it of filtration $\geqslant n$, so that it gives a cohomological invariant $a$ of $G$ of degree $n$, and one checks that $\< a,\sigma \> - \< a,\sigma' \> =1$; that information is enough to conclude the proof.
            
            \medskip
            
            \n {\it Dependence of} $\Inv$$_{k_0}(G)$ {\it on $H^\bullet(k_0)$ - Universal objects.}
            
            \smallskip
            \n (i) {\it Additive structure}
            
            \smallskip
            
            For the additive structure, $\Inv_{\C}$$(G)$ is a universal object, i.e., there is natural isomorphism of $\F$$_2$-vector spaces :
            $ {\rm Inv}$$_{k_0}(G)  \simeq {\rm Inv}_{\C}(G) \ \otimes_{{\bf F}_2} H^\bullet(k_0).$
            
            \smallskip
           
            \n (ii) {\it Ring structure}
            
            \smallskip
            
            For the ring structure, it is $\Inv_{\R}$$(G)$ which is a universal object : there is a natural graded-$\F$$_2$-algebra isomorphism  :
             $ {\rm Inv}_{k_0}(G)  \simeq \Inv_{\R}$$(G) \ \otimes_{H^\bullet(\R)} H^\bullet(k_0).$
             
             \smallskip
            
       \n     {\small[In this formula, $H^\bullet(k_0)$ is viewed as an $H^\bullet(\R)$-algebra via the unique homomorphism
            $H^\bullet(\R)$$ \to H^\bullet(k_0)$ which maps the class of $-1$ in $H^1(\R) \simeq \R^\times/(\R^\times)^2$ onto
 the class of $-1$ in $H^1(k_0) \simeq k_0^\times/(k_0^\times)^2$.] \par}
 
 \bigskip
 \begin{center}  References
 
 \end{center}

\n [B] N. Bourbaki, {\it Groupes et Alg\`{e}bres de Lie}, chap.IV-VI, Hermann, Paris, 1968.

\n [S] J-P. Serre, {\it Cohomological invariants, Witt invariants, and trace forms}, ULS 28, 1-100, AMS, 2003.

       \bigskip 
       
  \n {\it Note.} After my lecture, Stefan Gille has pointed out to me that, using a different method (based on a theorem of Totaro, but not involving involutions), Christian Hirsch had already computed in 2009 the structure of the cohomological invariants of all the finite Coxeter groups, under some mild hypotheses
  on the ground field; his method also applies to other types of invariants. Reference :
  
  \smallskip
  Christian Hirsch, {\it Cohomological invariants of reflection groups}, Diplomarbeit (Betreuer: Prof. Dr. Fabien Morel), Univ. München, 2009; available on arXiv\!:1805.04670[math.AG].
  
  \smallskip
\
       
       \bigskip
       
       \bigskip
       
       Jean-Pierre Serre
       
      Coll\`{e}ge de France
      
       3 rue d'Ulm
       
       75005 PARIS, France
   \end{document}